\documentclass[12pt]{article}

\newtheorem{theorem}{Theorem}[section]
\newtheorem{lemma}[theorem]{Lemma}
\newtheorem{proposition}[theorem]{Proposition}
\newtheorem{corollary}[theorem]{Corollary}

\newfam\msbfam
\font\tenmsb=msbm10  scaled \magstep1 \textfont\msbfam=\tenmsb
\font\sevenmsb=msbm7 scaled \magstep1 \scriptfont\msbfam=\sevenmsb
\font\fivemsb=msbm5  scaled \magstep1 \scriptscriptfont\msbfam=\fivemsb
\def\Bbb{\fam\msbfam \tenmsb}

\def\RR{{\Bbb R}}
\def\CC{{\Bbb C}}

\def\ss{\subseteq}
\def\ra{\rightarrow}

 \def\HollowBox #1#2{{\dimen0=#1 \advance\dimen0 by -#2       
       \dimen1=#1 \advance\dimen1 by #2                       
        \vrule height #1 depth #2 width #2                    
        \vrule height 0pt depth #2 width #1                   
        \llap{\vrule height #1 depth -\dimen0 width \dimen1}% 
       \hskip -#2                                             
       \vrule height #1 depth #2 width #2}}                   
 \def\BoxOpTwo{\mathord{\HollowBox{6pt}{.4pt}}\;}             
\def\endpf{\hfill $\BoxOpTwo$ \smallskip \\ }

\def\Re{\hbox{Re}\,}

\def\dbar{\overline{\partial}}

\newfam\msbbfam
\font\tenmsbb=msbm10  scaled \magstep1 \textfont\msbbfam=\tenmsbb
\font\sevenmsbb=msbm7  scaled \magstep1 \scriptfont\msbbfam=\sevenmsbb
\font\fivemsbb=msbm5    scaled \magstep1 \scriptscriptfont\msbbfam=\fivemsbb

\usepackage{graphicx}

\usepackage{amsmath}

\begin{document}

\begin{center}
\Large \bf Smoothness to the Boundary of Biholomorphic Mappings\footnote{{\bf Subject 
Classification Numbers:}  32H40, 32H02.}
\end{center}
\vspace*{.12in}

\begin{center}
Steven G. Krantz
\end{center}

\date{\today}

\begin{quote}
{\bf Abstract:}  Under a plausible geometric hypothesis, we show that a biholomorphic mapping
of smoothly bounded, pseudoconvex domains extends to a diffeomorphism of the closures.
\end{quote}
\vspace*{.25in}

\markboth{STEVEN G. KRANTZ}{SMOOTHNESS OF BIHOLOMORPHIC MAPPINGS}

\section{Introduction}

The Riemann mapping theorem (see [GRK]) tells us that the
function theory of a simply connected, planar domain $\Omega$,
other than than the entire plane, can be transferred from
$\Omega$ to the unit disc $D$.	But, for many questions, one
needs to know the behavior of the Riemann mapping at the boundary.

The first person to take up this issue was P. Painlev\'{e}. He
proved in his thesis that, if the domain $\Omega$ has $C^\infty$ boundary,
then the Riemann mapping (and its inverse) extends smoothly to
the boundary (see [BUR] for details of this history).  Later O. Kellogg 
gave a proof of this result that connected the Riemann mapping with 
potential theory.   Further on, Stefan Warschawski refined Kellogg's results
and gave substantive local boundary analyses of the Riemann mapping.

It was quite some time before any progress was made on this question
in the context of several complex variables.  The first real theorem
of a general nature was proved by C. Fefferman [FEF].   He showed
that a biholomorphic mapping of smoothly bounded, strongly pseudoconvex
domains in $\CC^n$ extends to a diffeomorphism of the closures.
Fefferman's work opened up a flood of developments in this subject.
We only mention here that Bell [BEL1] and Bell/Ligocka [BELL] were
able to greatly simplify Fefferman's proof by connecting the
problem in a rather direct fashion with the Bergman projection.
The work of Bell and Bell/Ligocka led to a number of simplifications, generalizations, and
extensions of Fefferman's result.   Many different mathematicians
have contributed to the development of this work.

The big remaining open problem is this:
\begin{quote}
{\bf Problem:}  Let $\Omega_1$ and $\Omega_2$ in $\CC^n$ be smoothly bounded, (weakly Levi) pseudoconvex
domains.  Let $\Phi: \Omega_1 \ra \Omega_2$ be a biholomorphic mapping.
Does $\Phi$ extend to a diffeomorphism of the closures?
\end{quote}
There are some counterexamples to this question---see for instance [FRI]---but
these definitely do not have smooth boundary.  In fact they do not even
have $C^2$ boundary.

In the present paper we are unable to give a full answer to this main problem.  But
we present the following somewhat encouraging partial result.

\begin{theorem} \sl
Let $\Omega_1$, $\Omega_2$ be smoothly bounded, Levi pseudoconvex domains in $\CC^n$.
Let $\Phi: \Omega_1 \ra \Omega_2$ be a biholomorphic mapping.  Assume that
$\Phi$ and $\Phi^{-1}$ each satisfy a Lipschitz condition of order exceeding
$(n-1)/n$.  Then $\Phi$ continues to a diffeomorphism of the closures of the domains.
\end{theorem}

\begin{corollary} \sl
Let $\Omega_1$, $\Omega_2$ be smoothly bounded, pseudoconvex domains in $\CC^n$.
Let $\Phi: \Omega_1 \ra \Omega_2$ be a biholomorphic mapping.  Assume that
$\Phi$ and $\Phi^{-1}$ each satisfy a Lipschitz condition of order 
1.  Then $\Phi$ continues to a diffeomorphism of the closures of the domains.
\end{corollary}

This result is in the nature of a bootstrapping result from partial differential
equations.  It seems to be the first general result---for {\it all} pseudoconvex domains---of
its kind.  And it has some basis in the history of the subject.  For Painlev\'{e} 
proved his theorem in dimension one by first establishing a result for $C^1$ boundary smoothness
of the mapping, and then bootstrapping.
No less an {\it eminence gris} than Jacques Hadamard cast public doubt on Painlev\'{e}'s
bootstrapping argument, and Painlev\'{e} had to work strenuously to defend his theorem.  
See [BUR] for the full history.

It may be noted that the hypothesis of Lipschitz continuity in the theorem is a nontrivial
one.  Henkin [HEN] was able to show, prior to Fefferman's celebrated result,
that a biholomorphic mapping of smoothly bounded, strongly pseudoconvex domains will extend
to be Lipschitz $1/2$ to the boundary.  He did so by analyzing and estimating
the Carath\'{e}odory metric.  But there are not many results of this type.

We see that the Lipschitz condition in Theorem 1.1 in case $n=2$ meshes
nicely with Henkin's result described in the last paragraph.

A final, rather significant, comment is this.  Our arguments here 
are inspired by those in [BEL1].  Bell uses global regularity ideas of
Kohn which exploit weighed $L^2$ spaces.  In the paper [BEL1],
a good deal of the work is expended in proving that the complex Jacobian
determinant $u$ of the mapping $\Phi$ is bounded.  This fact is used in
turn to prove that the complex Jacobian determinant $U$ of the inverse
mapping $\Phi^{-1}$ is nonvanishing.  As we shall see below, our hypothesis
of Lipschitz continuity of order exceeding $(n-1)/n$ obviates these arguments
and gets to the necessary result rather quickly.  The remaining steps 
comprise just one paragraph on page 108 of [BEL1].  We have to work
a bit harder because we need to set things up in the context of Kohn's
weighted $L^2$ spaces.  But the spirit of our arguments follows Bell.

We also warn the reader of the following point. The main thrust
of this paper is to prove estimates on the derivatives of the
mappings $\Phi$ and $\Phi^{-1}$. However our crucial Lemma
4.2, based on an old idea of S. R. Bell, entails taking a good
many derivatives of $\Phi$. So it appears as though there are
a number of extraneous terms in our calculations that involve
derivatives of $\Phi$. But we will go to quite a lot of extra
trouble to find a means of absorbing those extra derivatives.
In the end they will all be accounted for, and we will obtain
valid estimates for the derivatives of $\Phi$.

\section{Condition \boldmath $R$ and Related Ideas}

Throughout this paper we shall use the language of Sobolev spaces
(see [ADA], [STE]).  For $s$ a nonnegative integer and $1 \leq p \leq \infty$, 
we let $W_{s, p}$ denote the usual
Sobolev space of functions with $s$ weak derivatives in $L^p$.	The norm
that we use on the Sobolev space is standard, and we refer to [ADA] 
for details.

One of the important innovations that S. R. Bell introduced into this
subject is his {\it Condition $R$}.  It says this:
\begin{quote}
{\bf Condition \boldmath $R$:}	Let $\Omega \ss \CC^n$ be a smoothly
bounded domain.  We say that $\Omega$ satisfies {\it Condition  R} if
the Bergman projection $P$ maps $C^\infty(\overline{\Omega})$ to
$C^\infty(\overline{\Omega})$.  Equivalently, for each $s > 0$, there
is an $m(s) > 0$ so that the Bergman projection $P$ maps the Sobolev
space $W_{m(s), 2}(\Omega)$ to the Sobolev space $W_{s, 2}(\Omega)$.
\end{quote}

In what follows we shall suppose that $1 < m(1) < m(2) < \cdots \ra \infty$
and that each $m(j)$ is an integer.

In the paper [BEL1], Bell proves the following elegant result:
\begin{theorem} \sl
Let $\Omega_1$, $\Omega_2$ be smoothly bounded, pseudoconvex domains
in $\CC^n$.  Assume that one (but not necessarily both) of these domains
satisfies Condition $R$.  Then any biholomorphic mapping $\Phi: \Omega_1 \ra \Omega_2$
extends to a diffeomorphism of the closures.
\end{theorem}

\section{Ideas of Kohn}

The classical treatment of the $\overline{\partial}$-Neumann problem
is based on the traditional Euclidean $L^2$ inner product---see [FOK].  Kohn's
idea in [KOH] is to use an inner product with a weight.  This is
inspired by work of H\"{o}rmander [HOR], and that in turn comes from
old ideas of Carleman.

Kohn's setup is this (see [KOH, p. 279]).  Fix a smoothly bounded domain $\Omega$
in $\CC^n$.  Let $\lambda$ be a $C^\infty$, nonnegative
function on a neighborhood of $\overline{\Omega}$.  Usually $\lambda$
will be strictly plurisubharmonic.  With $\lambda$ fixed and $t \geq 0$,
we shall define the $\overline{\partial}$-Neumann problem of weight $t$, with
real $t > 0$.  We let ${\cal A}$ be the space of all forms on $\overline{\Omega}$ which
have $C^\infty$ coefficients up to the boundary.  For $\phi, \psi \in {\cal A}$, we define
$$
\langle \phi, \psi\rangle_{(t)} \equiv \langle \phi, e^{-t\lambda} \psi \rangle \qquad \hbox{and}
   \qquad \|\phi\|_{(t)}^2 = \langle \phi, \phi \rangle_{(t)} \, .
$$
Here $\langle \ , \ \rangle = \langle \, , \, \rangle_{(0)}$ is the usual $L^2$ inner product.

It is an easily verified fact that the norms $\| \ \ \|_{(t)}$ are equivalent to the norm
$\| \ \ \|_0 = \| \ \ \|$.  Hence a function is in the completion under any of these
norms if and only if it is square integrable.  We let $\widetilde{\cal A}_t$ be the Hilbert
space obtained by completing ${\cal A}$ under the norm $\| \ \ \|_{(t)}$.

The $\overline{\partial}$-Neumann problem may be set up in the $\langle \ , \ \rangle_{(t)}$ inner
product rather than the usual $L^2$ inner product $\langle \ , \ \rangle$.  These are familiar
ideas, and the details are provided in [KOH].  One of the main points that must be noted
is that the {\it formal adjoint} of the operator $\overline{\partial}$, when calculated
in the $\langle \ , \ \rangle_{(t)}$ inner product, is
$$
{\cal I}_t = {\cal I} - t \sigma ({\cal I}, d\lambda)\, .
$$
Here $\sigma$ is the ``symbol'' in the usual sense of pseudodifferential operators.
Also ${\cal I}$ is the standard formal adjoint of $\overline{\partial}$ with respect
to the standard Euclidean Hermitian inner product and ${\cal I}_t$ is the formal
adjoint of $\overline{\partial}$ with respect to the inner product $\langle \ , \ \rangle_{(t)}$.
We thus see how the parameter $t$ comes into play.  If $t$ is chosen positive and large
enough, then certain terms in the usual $\overline{\partial}$-Neumann estimates can
be forced to dominate certain others.  Again see [KOH] for the details.

Let $\Omega_1$, $\Omega_2$ be smoothly bounded, pseudoconvex domains and $\Phi: \Omega_1 \ra \Omega_2$
a biholomorphic mapping which is bi-Lipschitz of order exceeding $(n-1)/n$.  We shall apply the preceding ideas
on $\Omega_2$ with $\lambda(z) = \lambda_2(z) = |z|^2$ and on $\Omega_1$ with $\lambda(z) = \lambda_1(z) = |\Phi(z)|^2$.
Note that we are applying Kohn's construction twice.\footnote{It is because the weight $e^{- t |\Phi(z)|^2}$ gets
differentiated in the proofs below that we must be careful to absorb these error terms.}

In this context we shall refer to the $t$-weighted Bergman projection as $P_{t, 1}$ (for $\Omega_1$) and $P_{t, 2}$ (for $\Omega_2$).
We shall also call Bell's regularity condition in the context of Kohn's weighted inner
product by the name ``Condition $R_t$.''  We shall denote the Bergman kernels by $K_{t,1}$
and $K_{t,2}$  .
As a result of these ideas, the $\overline{\partial}$-Neumann problem on $\Omega_2$, formulated with the indicated
weight $\lambda$, satisfies favorable estimates (this follows from [KOH]) as long as $t$ is large enough.  
Hence $\Omega_2$ satisfies Condition $R_t$.  Bell also makes use of these facts.

These are the tools that we shall need in the next section to get to our result.

In what follows we shall take it that we are working with the Bergman theory for the
inner product $\langle \ , \ \rangle_{(t)}$ for $t$ sufficiently large, and that 
$\Omega_2$ satisfies Condition $R_t$.  We formulate this last by saying that $P_{t,2}: W_{m(s),2} (\Omega) \ra W_{s,2}(\Omega)$
for any $s \geq 0$ and suitable $m(s) \geq s$.

Sometimes, in what follows, we will talk about 
\smallskip \\

{\bf (i)} a domain $\Omega$ with weight $\lambda$ 
\smallskip \\

\noindent but make no reference to
\smallskip \\

{\bf (ii)} $\Omega_1$, $\Omega_2$, or the mapping $\Phi$.  
\smallskip \\

\noindent We will later apply {\bf (i)} to {\bf (ii)}.

\section{The Guts of the Proof}

In this section, $\Omega$ is a smoothly bounded, pseudoconvex domain equipped with
the weight $e^{-t \lambda}$.

\begin{lemma} \sl
Let $\Omega \subset \! \subset \CC^n$ be smoothly bounded and pseudoconvex.  Suppose that the $\lambda$ from
the weight on $\Omega$ is smooth on $\Omega$.  Assume that
$\Omega$ satisfies Condition $R_t$ with respect to the projection $P_t$.  Let $w \in \Omega$ be
fixed.  Let $K_t$ denote the Bergman kernel.   Then there is a constant $C_w > 0$ so that
$$
\|K_t(w, \, \cdot \, )\|_{\rm sup} \leq C_w \, .
$$
\end{lemma}
\noindent {\bf Proof:}  The function $K(z, \, \cdot \, )$ is harmonic.  Let $\varphi: \Omega \ra \RR$ 
be a radial, $C^\infty$ function centered at $w$ and supported in $\Omega$ so that the
radius of the support is comparable to half the distance of $w$ to the boundary.  Assume that $\varphi \geq 0$ and $\int \varphi(\zeta) \, dV(\zeta) = 1$.
Then the mean value property implies that
$$
K_t(z, w) = \int_\Omega K_t(z, \zeta) \varphi(\zeta) \, dV(\zeta) 
            = \int_\Omega K_t(z, \zeta) \bigl [ \varphi(\zeta) e^{t\lambda(\zeta)} \bigr ] e^{-t\lambda(\zeta)} \, dV(\zeta)  \, .
$$
Of course this last displayed expression equals $P_t \biggl (\varphi(\, \cdot \, ) e^{t\lambda(\, \cdot \, )} \biggr )$.  Thus
\begin{eqnarray*}
\|K_t(w, \, \cdot \, ) \|_{\rm sup} & = & \sup_{z \in \Omega} |K_t(w, z)| \\
                                  & = & \sup_{z \in \Omega} |K_t(z, w)| \\
				  & = & \sup_{z \in \Omega} |P_t \bigl [\varphi(\, \cdot \, ) e^{t\lambda(\, \cdot \, )} \bigr ]| \, .
\end{eqnarray*}

By Sobolev's theorem, this last is 
$$
\leq C(\Omega, w) \|P_t \bigl [\varphi(\, \cdot \, ) e^{t\lambda(\, \cdot \, )} \bigr ]\|_{W_{2n+1,2}} \, .
$$
By Condition $R_t$, this is 
$$
\leq C(\Omega, w) \cdot \|\varphi(\, \cdot \, ) e^{t\lambda(\, \cdot \, )} \|_{W_{m(2n+1),2}} \equiv C_w \, .   \eqno \BoxOpTwo
$$

\noindent {\bf Remark:} It is worth noting that the estimate
obtained in this last proof depends on some derivatives of
$\lambda$ on a compact set. In practice this causes no harm.
We only need to know that $\|K_t(w, \, \cdot \, )\|_{\rm sup}$
is bounded so that we can perform an integration by parts in
the proof of Lemma 4.2 below. 
\smallskip \\

\begin{lemma} \sl
Let $u \in C^\infty(\overline{\Omega})$ be arbitrary.  Let $s \in \{0, 1, 2, \dots\}$.  Then
there is a $v \in C^\infty(\overline{\Omega})$ such that $P_t v = 0$ and the functions $u$ and
$v$ agree to order $s$ on $\partial \Omega$.
\end{lemma}
\vspace*{.15in}

\noindent {\bf Remark:}  This lemma in the present formulation is not entirely
satisfactory.   For, in its statement here, we suppose that the
weight $\lambda$ is smooth across the boundary.  Yet, in the applications below,
the weight is taken to be $|\Phi(z)|^2$, and that is {\it not} known a priori
to be smooth up to the boundary (in fact our goal is to prove that it is smooth
up to the boundary).  

The way to address this problem is to use an approximation argument.  In the case
that the domain $\Omega_1$ is convex, the approximation is very simple.  We simply
replace $\Omega_1$ by $(1 - \epsilon)\Omega_1$, $\epsilon > 0$, so that the mapping is
smooth across the boundary.  The relevant estimates are uniform in $\epsilon$, and the
result is correct in the limit.

For non-convex $\Omega_1$, we must take advangage of ideas in [BEL2].  For Bell explains
there how to localize the smoothness-to-the-boundary arguments that we present here.  As 
a result, if $p \in \partial \Omega_1$, $\nu_p$ is the Euclidean unit outward normal vector
at $p$, and $U$ is a small neighborhood of $p$, then we may apply our arguments on $\Omega_1^\epsilon \equiv (U \cap \Omega_1) - \epsilon \nu_p$.
Thus the mapping will be smooth across the boundary and (a localized version of) Lemma 4.2 will apply without
any problem.  All the relevant estimates are uniform in $\epsilon$, and our desired result
holds in the limit.
\medskip \\

\noindent {\bf Proof:}   This lemma is the key to Bell's approach to these matters.  
We will need to expend some effort to adapt Bell's ideas to the new weighted context.

We of course assume that our domain $\Omega$ is equipped with an inner
product $\langle \ , \ \rangle_{(t)}$ based on a weight $e^{- t \lambda}$.

After a partition of unity, it suffices to prove the assertion
in a small neighborhood $W$ of $z_0 \in \partial \Omega$.  After a rotation,
we may assume that $\partial \rho/\partial z_1 \ne 0$
on $W \cap \overline{\Omega}$.  [Here $\rho$ is a defining function
for the domain $\Omega$---see [KRA1].]

Define the differential operator
$$
\nu = \frac{\Re \left \{ \sum_{j=1}^n \frac{\partial \rho}{\partial z_j} \frac{\partial}{\partial z_j} \right \} }{\sum_{j=1}^n \left | \frac{\partial \rho}{\partial z_j} \right |^2} \, .
$$
Observe that $\nu \rho \equiv 1$.  Now we shall define $v$ by induction on $s$.
In what follows, we shall make use of the differential operator
$$
T = \frac{\partial}{\partial \zeta_1} + t \frac{\partial \lambda}{\partial \zeta_1} \, .
$$

For the case $s = 0$, we set
$$
w_1 = \frac{\rho u}{T \rho} \, .
$$
If $W$ is small then of course $T \rho$ will not vanish.
Also define
$$
v_1 = T w_1 = u + {\cal O}(\rho) \, .
$$

Then we see immediately that $u$ and $v_1$ agree to order 0 on $\partial \Omega$.  Furthermore,
\begin{eqnarray*}
P_t v_1 & = & \int K_t(z, \zeta) T w_1 e^{-t\lambda} \, dV \\
	& = & - \int T_\zeta \left [ K_t(z,\zeta) e^{-t \lambda} \right ] w_1 \, dV \\
	& = & 0 \, .
\end{eqnarray*}
The penultimate equality comes from integration by parts.   
This operation is justified by Lemma 4.1.  Note that $T_\zeta$
annihilates $K_t(z,\zeta) e^{-t\lambda(\zeta)}$ by a simple calculation (using
the fact that $K_t(z,\zeta)$ is conjugate holomorphic in the $\zeta$ variable).

Now suppose inductively that 
\begin{eqnarray*}
w_{s-1} & = & w_{s-2} + \theta_{s-1} \cdot \rho^{s-1} \, , \\
v_{s-1} & = & T w_{s-1} 
\end{eqnarray*}
(for some smooth function $\theta_{s-1}$).
We construct 
$$
w_s = w_{s-1} + \theta_s \cdot \rho^s
$$
so that
$$
v_s \equiv T w_s
$$
agrees to order $s-1$ with $u$ on $\partial \Omega$.

By the inductive hypothesis,
\begin{eqnarray*}
v_s & = & T w_s \\
    & = & T w_{s-1} + T(\theta_s \rho^s) \\
    & = & v_{s-1} + \rho^{s-1} \left [ s \theta_s T \rho + \rho T \theta_s \right ] \, .
\end{eqnarray*}
This expression agrees, by the inductive hypothesis, with $u$ to order $s-1$ on $\partial \Omega$.  We now
must examine $D (u - v_s)$, where $D$ is any $s$-order differential operator.  There are two cases:
\begin{description}
\item[{\bf Case 1:}]  Assume that $D$ involves a tangential derivative $D_0$.  Then we may write
$D = D_0 D_1$.  Then
$$
D (u - v_s) = D_0 \alpha \, ,
$$
where $\alpha$ vanishes on $\partial \Omega$.  But then it follows that $D_0 \alpha = 0$ because
$D_0$ is tangential.  
\item[{\bf Case 2:}]  Now assume that $D$ has no tangential derivative in it.  So we take
$D = \nu^s$, where $\nu$ was defined at the beginning of this discussion.  Our job is to
choose $\theta_s$ so that
$$
\nu^s (u - v_s) = 0 \quad \hbox{on} \ \partial \Omega \, .
$$
So we require that
$$
\nu^s(u - v_{s-1}) - \nu^s \left ( T(\theta_s \rho^s) \right ) = 0 \quad \hbox{on} \ \partial \Omega \, .
$$
This is the same as
$$
\nu^s (u - v_{s-1}) - \theta_s \left ( \nu^s T \rho^s \right ) = 0 \quad \hbox{on} \ \partial \Omega
$$
(because terms that contain $\rho$ must vanish on $\partial \Omega$) or
$$
\nu^s(u - v_{s-1}) - \theta_s \left (\nu^s \frac{\partial}{\partial \zeta_1} \rho^s \right )
              - \theta_s \left ( \nu^s t \frac{\partial \lambda}{\partial \zeta_1} \rho^s \right ) = 0 \quad \hbox{on} \ \partial \Omega \, .
$$
This may be rewritten as
$$
\nu^s (u - v_{s-1}) - \theta_s \cdot s! \frac{\partial \rho}{\partial \zeta_1} 
            - t \cdot \theta_s \cdot \tau \cdot \rho  \, ,
$$
where $\tau$ stands for terms that come from the differentiations.
The last line may be rewritten as
$$
\theta_s = \frac{\nu^s(u - v_{s-1})}{s! \frac{\partial \rho}{\partial \zeta_1} + t \cdot \tau \cdot \rho} \, .
$$
If $W$ is small enough then the denominator cannot vanish and we see that we have a well-defined choice for $\theta_s$
as desired.
\endpf
\end{description}

We note that, in [BEL1], Bell has a particularly elegant way of expressing the content of this last lemma.
His formulation will be useful for us later, so we formulate it now.  First some notation.

If $\Omega \ss \CC^n$ is a domain (a connected, open set), then let $W_{s,p}^0(\Omega)$ be the closure
of $C_c^\infty(\Omega)$ in $W_{s,p}(\Omega)$.   Now we have Bell's formulation of our Lemma 4.2:

\begin{corollary} \sl
Let $\Omega$ be a smoothly bounded, pseudoconvex domain.  Let $s,m \in \{0,1,2,\dots\}$.  Then there is a linear
operator $\Psi^{s,m}: W_{s+m, 2}(\Omega) \ra W_{s,2}^0(\Omega)$ such that $P_t \Psi^{s,m} = \hbox{id}$.
The norm of this operator depends polynomially on $t$ and on derivatives of $\lambda$.
\end{corollary}			     

\section{A Deeper Analysis}

A troublesome feature of Lemma 4.2 and Corollary 4.3 is that the weight $\lambda$ gets differentiated $s$ times,
and $\lambda$ (in practice) is defined in terms of the mapping $\Phi$.  Since our job in the end is to estimate the derivatives
of $\Phi$, this looks problematic.  We need to develop a way to absorb these extraneous derivatives of $\Phi$.

With that thought in mind, we recall the standard Sobolev embedding theorem for a smooth domain in $\RR^N$ (see [ADA],
[STE] for details).

\begin{proposition}[Sobolev] \sl
Let $\Omega \ss \RR^N$ be a smoothly bounded domain.   Let $W_{m,p}$ be the standard
Sobolev space of functions on $\Omega$ having $m$ weak derivatives in the space $L^p$.
Equip $W_{m,p}$ with the usual norm.  Then we have the embedding
$$
W_{k,p} \ss W_{\ell, q}
$$
whenever $k > \ell$, $1 \leq p < q \leq \infty$, and
$$
\frac{1}{q} = \frac{1}{p} - \frac{k - \ell}{N} \, .
$$
\end{proposition}

We are particularly interested in domains in $\CC^n$.  Hence, for us, $N = 2n$.  Also
we will apply the result in case $k = 4$ and $\ell = 2$.  Thus the
important inclusion is
$$
W_{k + n/2, 2} \ss W_{k, 4}\, .
$$
We will generally exploit this embedding in the form of the inequality
$$
\| f \|_{k, 4} \leq C \|f\|_{k + n/2, 2} \, .
$$

We will also make good use of the following refinement of the Sobolev theorem
that is due to Ehrling, Gagliardo, and Nirenberg (see [ADA, p.\ 75] for the details):

\begin{theorem} \sl
Let $\Omega \ss \RR^N$ be a smoothly bounded domain.   Let $W_{m,p}$ be the standard
Sobolev space of functions on $\Omega$ having $m$ weak derivatives in the space $L^p$.
Equip $W_{m,p}$ with the usual norm.  Let $\epsilon_0 > 0$.  Let $m$ be a positive integer.
Then there is a constant
$K$, depending on $\epsilon_0$, $m$, $p$, and $\Omega$ so that, for any integer
$j$ with $0 \leq j \leq m-1$, any $0 < \epsilon < \epsilon_0$, and any $u \in W_{m, p}$,
$$
\|u\|_{j,p} \leq K \epsilon \|u\|_{m, p} + K \epsilon^{-j/(m-j)} \|u\|_{0,p} \, .
$$
\end{theorem}

It is worth taking some time now to do a little analysis.  Examine the proof 
of Lemma 4.2. [Note that, when we apply Lemma 4.2 and Corollary 4.3, we do so on $\Omega_1$
with $\lambda(z) = \lambda_1(z) = |\Phi(z)|^2$.]  At each stage we integrate by parts, and therefore a derivative
falls on $\lambda$ (and hence on $\Phi$).  Thus we see that the function $v$ that
we construct has the form
$$
v = q_0 + q_1 t\nabla \Phi + q_2 t^2\nabla^2 \Phi + \cdots + q_s t^s \nabla^s \Phi
        + q'_1 t\nabla \overline{\Phi} + q'_2 t^2\nabla^2 \overline{\Phi} + \cdots + q'_s t^s \nabla^s \overline{\Phi}    \, .
$$
Here we use the notation $\nabla \Phi$ or $\nabla \overline{\Phi}$ to denote some derivative of some component of $\Phi$ 
or $\overline{\Phi}$ and
$\nabla^j \Phi$ or $\nabla^j \overline{\Phi}$ to denote some $j^{\rm th}$ derivative of some component of $\Phi$ or $\overline{\Phi}$.  
Also $q_j$, $q'_j$ denotes an expression that involves components of $\Phi$, derivatives of $\rho$, and
derivatives of $\theta_j$.
Thus we see that
\begin{eqnarray*}
  \|v\|_{r, 2} & \leq & C \cdot \bigl [ \int \bigl |q_0 \bigr |^2 \, dV^{1/2} + \int \bigl | \nabla^r q_0 \bigr |^2 \, dV^{1/2} 
          + \int |q_1 t \nabla\Phi|^2\, dV^{1/2}  \\
               &      & \quad + \int \bigl |(t \nabla^r q_1) \Phi \bigr |^2 + \bigl | t q_1 \nabla^{r+1} \Phi \bigr |^2 \, dV^{1/2} \\
               &      & \quad + \int |q_2 t \nabla^2\Phi|^2\, dV^{1/2}  
 + \int \bigl | (t^2 \nabla^r  q_2) \Phi \bigr |^2 + \bigl | t^2 q_2 \nabla^{r+2} \Phi \bigr |^2 \, dV^{1/2} \\
               &      & \qquad + \cdots + \int |q_s t \nabla^s \Phi|^2\, dV^{1/2}  \\
               &      & \qquad + \int \bigl | (t^s \nabla^r  q_s) \Phi \bigr |^2 + \bigl | t^s q_s \nabla^{r+s} \Phi \bigr |^2 \, dV^{1/2} \bigr ] \\
               &      & \ \ \hbox{plus similar terms involves $q'_j$ and $\nabla^j \overline{\Phi}$.} \\
\end{eqnarray*}		

Using H\"{o}lder's inequality, we see that this last is majorized by
\begin{eqnarray*}
\lefteqn{C \cdot \biggl [ \int |q_0|^2 \, dV^{1/2} + \int \bigl | \nabla^r q_0\bigr |^2 dV^{1/2}} \\
                          && + \quad t \int |q_1|^4\, dV^{1/4} \cdot \int|\nabla \Phi|^4\, dV^{1/4}  + t \int \bigl | \nabla^r q_1 \bigr |^4 \, dV^{1/4} \cdot \int \bigl | \Phi \bigr |^4 \, dV^{1/4}  \\
                          && + t\int \bigl | q_1 \bigr |^4 \, dV^{1/4} \cdot \int \bigl | \nabla^{r+1} \Phi \bigr |^4 \, dV^{1/4}  \\
		   && \quad + \quad t^2 \int |q_2|^4\, dV^{1/4} \cdot \int|\nabla^2 \Phi|^4\, dV^{1/4} + t^2 \int \bigl | \nabla^r q_2 \bigr |^4 \, dV^{1/4} \cdot \int \bigl | \Phi \bigr |^4 \, dV^{1/4}  \\
                   && \qquad + t^2 \int \bigl | q_2 \bigr |^4 \, dV^{1/4} \cdot \int \bigl | \nabla^{r+2} \Phi \bigr |^4 \, dV^{1/4} + \cdots \\
		   && \qquad + \quad t \int |q_s|^4\, dV^{1/4} \cdot \int|\nabla^s \Phi|^4\, dV^{1/4} + t^s \int \bigl | \nabla^r q_s \bigr |^4 \, dV^{1/4} \cdot \int \bigl | \Phi \bigr |^4 \, dV^{1/4} \\
                   && \qquad + t^s \int \bigl | q_s \bigr |^4 \, dV^{1/4} \cdot \int \bigl | \nabla^{r+s} \Phi \bigr |^4 \, dV^{1/4} \biggr ] \\
		   &      & \ \ \hbox{plus similar terms involves $q'_j$ and $\nabla^j \overline{\Phi}$.} \\
\end{eqnarray*}
Now we may estimate this more elegantly as
$$
C' \cdot  (1 + t)^s \bigl [ 1 + \|\Phi\|_{r+s,4} \bigr ] \, .
$$
And then we apply the Sobolev inequality noted above to estimate this at last by
$$
C' \cdot  (1 + t)^s \bigl [1 + \|\Phi\|_{r + s + n/2, 2} \bigr ] \, .
$$
At last we apply the Ehrling-Gagliardo-Nirenberg inequality (with $m = r + s + n$, $j = r + s + n/2$) to obtain
$$
\leq C''(1 + t)^s \bigl (1 + \epsilon \cdot \|\Phi\|_{r + s + n, 2} + \epsilon^{-(r + s + n/2)/(n/2)} \cdot \|\Phi\|_{0,2} \bigr ) \, .     \eqno (5.2)
$$
This inequality is valid for any $0 < \epsilon < 1$ provided that $C''$ is large enough.   Furthermore,
$C''$ depends on $s$.

We finally note that, in the proof of Lemma 4.2, the definition of $w_1$ involves a division by $T \rho$
(and the definition of $T$ entails a derivative of $\Phi$).  But this can be treated with a Neumann series,
or just using the quotient rule.

Now we can reformulate Corollary 4.3 as follows:  

\setcounter{theorem}{2}

\begin{corollary} \rm
Let $\Omega$ be a smoothly bounded, pseudoconvex domain.  Let $s \in \{0,1,2,\dots\}$.  Fix
any $\epsilon > 0$.  Then there is a linear
operator $\Psi^s: W_{r + s + n,2}(\Omega) \ra W_{r,2}^0(\Omega$ such that $q_t \Psi^s = \hbox{id}$.  Moreover,
the norm of this operator does not exceed $C''(1 + t)^s \bigl (1 + \epsilon \cdot \|\Phi\|_{r + s + n, 2} + \epsilon^{-(r + s + n/2)/(n/2)} \cdot \|\Phi\|_{0,2} \bigr )$.
\end{corollary}

\section{The Final Argument}

For the rest of this discussion, we let $\Omega_1$ and $\Omega_2$ be fixed,
smoothly bounded, pseudoconvex domains in $\CC^n$. We fix a strictly
plurisubharmonic function $\lambda(z) = \lambda_2(z) = |z|^2$, and we equip $\Omega_2$ with the inner product with weight
$e^{-t \lambda(z)}$.  We also, as above, equip $\Omega_1$ with the inner
product having weight $\lambda(z) = \lambda_1(z) = |\Phi(z)|^2$.  We assume that there is a biholomorphic mapping
$\Phi: \Omega_1 \ra \Omega_2$ that extends in a bi-Lipschitzian fashion,
of order greater than $(n-1)/n$, to
the boundary, and we equip $\Omega_1$ with the inner product with weight
$e^{-t \lambda |\Phi(z)|^2}$. We let $P_{t,1}$ and $P_{t,2}$ be the resulting
Bergman projections for $\Omega_1$, $\Omega_2$ respectively. We let $u$
denote the complex Jacobian determinant of $\Phi$. And we let $U$ denote
the complex Jacobian determinant of $\Phi^{-1}$. So $u$ is a
complex-valued holomorphic function on $\Omega_1$ and $U$ is a
complex-valued holomorphic function on $\Omega_2$. Finally, for $j = 1,
2$, we let $\delta_j(z) = \delta_{\Omega_j}(z) = \hbox{dist}_{\rm Euclid} (z, {}^c
\Omega_j)$ for $z \in \Omega_j$.

\begin{lemma} \sl
The Bergman kernels for the two domains are related by
$$ 
K_{t,1}(z, \zeta) = u(z) \cdot K_{t,2}(\Phi(z), \Phi(\zeta)) \cdot \overline{u(\zeta)} \, .   \eqno (6.1.1)
$$
\end{lemma}
\noindent {\bf Proof:}  This is a standard change-of-variables argument, using the canonical
relationship between the real Jacobian determinant of a biholomorphic mapping
and the complex Jacobian determinant of that mapping (see Lemma 1.4.10 of [KRA1]).

Now, if $f$ is a Bergman space function on $\Omega_1$, then we have
\begin{eqnarray*}
\lefteqn{\int_{\Omega_1} \bigl [ u(z) K_{t,2}(\Phi(z), \Phi(\zeta)) \overline{u(\zeta)} \bigr ] f(\zeta) e^{-t\lambda(\Phi(\zeta)} \, dV(\zeta)}  \\
      & = & \int_{\Omega_2} u(z) K_{t,2}(\Phi(z), \xi) \overline{u(\Phi^{-1}(\xi))}  \\
      && \quad \times  f(\Phi^{-1}(\xi)) u^{-1}(\xi) \overline{u^{-1}(\xi)} e^{-t\lambda(\xi)} \, dV(\xi) \\
      & = & f(z) u(z) u^{-1}(\Phi(z)) \\
      & = & f(z) \, .
\end{eqnarray*}
Thus we see that the righthand side of (6.1.1) has the reproducing property on $\Omega_1$.  It is also conjugate symmetric and is
an element of the Bergman space in the first variable.  Therefore it must equal $K_{t,1}(z, \zeta)$.
\endpf \\

\begin{lemma} \sl
For any function $g \in L^2(\Omega_2)$, we have
$$
P_{t,1}\left ( u \cdot (g \circ \Phi) \right ) = u \cdot \left ( ( P_{t,2}(g) \circ \Phi \right ) \, .
$$
\end{lemma}
\noindent {\bf Proof:}  We use the preceding lemma to calculate that 
\begin{eqnarray*}
u(z) \cdot \left ( ( P_{t,2}(g) \circ \Phi \right ) (z) & = & 
                 u(z) \int_{\Omega_2} K_{t,2}(\Phi(z), \zeta) g(\zeta) e^{-t\lambda(\zeta)} \, dV(\zeta) \\
		 & = & u(z) \int_{\Omega_2} u(z)^{-1} K_{t,1}(z, \Phi^{-1}(\xi)) \overline{u(\Phi(\xi))^{-1}} \\
		 && \quad \times g(\zeta) e^{-t\lambda(\zeta)} \, dV(\zeta) \\
		 & = & u(z) \int_{\Omega_1} u(z)^{-1} K_{t, 1} (z, \xi) \overline{u(\xi)^{-1}}  \\
                 && \quad \times g(\Phi(\xi)) e^{-t\lambda(\Phi(\xi))} u(\xi) \overline{u(\xi)} \, dV(\xi) \\
		 & = & \int_{\Omega_1} K_{t,1}(z, \xi) g(\Phi(\xi)) u(\xi) e^{-t\lambda(\Phi(\xi))} \, dV(\xi) \\
		 & = & P_{t,1}\left ( u \cdot (g \circ \Phi) \right ) (z) \, .
\end{eqnarray*}
That establishes the result.
\endpf \\

It will be useful to have the following corollary, in which $\Phi$ is replaced by $\Phi^{-1}$ \
(and of course the corresponding Bergman kernels switch roles):

\begin{corollary} \sl
Let $U$ denote the complex Jacobian determinant of $\Phi^{-1}$.
Then, for any function $g \in L^2(\Omega_1)$, we have
$$
P_{t,2}\left ( U \cdot (g \circ \Phi^{-1}) \right ) = U \cdot \left ( ( P_{t,1}(g) \circ \Phi^{-1} \right ) \, .
$$
\end{corollary}

\begin{lemma} \sl
Let $H^\infty(\overline{\Omega}_1)$ denote the space of holomorphic functions on $\Omega_1$
which extend smoothly to $\overline{\Omega}_1$.   Let $s \in \{0,1,2,\dots\}$.  If $h \in H^\infty(\overline{\Omega}_1)$, then let $\phi_s = \Psi^s h$, where
$\Psi^s$ is introduced in Corollary 4.3.  Then
$$
U \cdot (h \circ \Phi^{-1}) = P_{t,2} (U \cdot (\phi_s \circ \Phi^{-1})) \, .
$$
\end{lemma}
\noindent {\bf Proof:}  We calculate, using Corollary 6.3, that
$$
P_{t,2}(U \cdot (\phi_s \circ \Phi^{-1})) = U \cdot (P_{t,1}(\phi_s) \circ \Phi^{-1}) =
          U \cdot (P_{t,1}(\Psi^sh) \circ \Phi^{-1}) = U \cdot (h \circ \Phi^{-1}) \, .  \eqno \BoxOpTwo
$$
\vspace*{.12in}

The next lemma has nothing to do with Condition $R$.  It is really only calculus.

\begin{lemma} \sl
Suppose that $\Phi^{-1}: \Omega_2 \ra \Omega_1$ is a biholomorphic mapping
between smoothly bounded, pseudoconvex domains in $\CC^n$.  Assume that $\Phi$ is bi-Lipschitz
of order exceeding $(n-1)/n$.  Let $U$ denote
the complex Jacobian determinant of $\Phi^{-1}$.  For each nonnegative integer $s$,
there is an integer $N = N(s)$ such that the operator
$$
g \longmapsto U \cdot (g \circ \Phi^{-1})
$$
is bounded from $W_{s + N,2}^0(\Omega_1)$ to $W_{s,2}^0(\Omega_2)$.
\end{lemma}
\noindent {\bf Proof:}  In what follows we let $d_j(z)$ denote the Euclidean
distance of $z$ from the boundary of $\Omega_j$.

Since the components of $\Phi^{-1}$ are holomorphic and Lipschitz
of order exceeding $(n-1)/n$, the derivatives of $\Phi^{-1}$ satisfy finite growth conditions at the boundary (see[GOL]).
That is to say
$$
\left | \frac{\partial^\alpha \Phi^{-1}}{\partial w^\alpha} (w) \right | \leq C \cdot d_2(w)^{-k + (n-1)/n} \, .  \eqno (6.5.1)
$$
Here $\alpha$ is a multi-index, $k = |\alpha|$, and $d_j$ is the distance of the argument
to the boundary of $\Omega_j$, $j = 1,2$.  Estimates like this one go back to Hardy and Littlewood
(see [GOL]).  

Now Sobolev's lemma and Taylor's formula tell us that, for $g \in W_{s+|\alpha| + n,2}^0(\Omega_1)$, 
$$
|D^\alpha g(z)| \leq C \cdot \|g\|_{s + |\alpha| + n} d_1(z)^s \, .
$$
For a given $s$, in order to show that
an $N$ exists so that $g \mapsto U \cdot (g \circ \Phi^{-1})$ is bounded from $W_{s+N,2}^0(\Omega_1)$
to $W_{s,2}^0(\Omega_2)$, it will suffice to show
that there is an integer $m > 0$ such that $d_1(\Phi^{-1}(w))^m \leq C \cdot d_2(w)$.
That such an $m$ exists is proved by Range [RAN1, Lemma 3.1].  The proof, naturally,
consists of applying Hopf's lemma to $\rho \circ \Phi^{-1}$, where $\rho$ is a bounded, plurisubharmonic
exhaustion function for $\Omega_2$ of the form $v d_2^{1/m}$ with $v \in C^\infty(\overline{\Omega}_2)$
and $v < 0$ on $\overline{\Omega}_2$.   Of course Diederich and Forn\ae ss [DIF]
have proved the existence of such an exhaustion function.  Range [RAN2] has given
a simpler approach to the matter, with the penalty of assuming greater boundary smoothness.

That completes the proof of the lemma.
\endpf
\smallskip \\

\begin{lemma} \sl
Let $s \in \{0,1,2,\dots\}$.  With notation as above,
$$
\|U \cdot (h \circ \Phi^{-1})\|_s \leq \|h\|_{s+N} \, .
$$
\end{lemma}
\noindent {\bf Proof:}  We note that Kohn's theory (see [BEL1] for the details)
entails that $P_{t,2}$ maps $W_{s, 2}(\Omega_2)$ to $W_{s,2}(\Omega_2)$ for $t$ sufficiently
large and any $s$.

Now we apply Corollary 6.3 and then Lemma 6.4 to see that
\begin{eqnarray*}
\|U \cdot (h \circ \Phi^{-1}) \|_s & \leq & \|P_{t,2} (U \cdot (\phi_s \circ \Phi^{-1} )) \|_s \\
                                   & \leq & \|U \cdot (\phi_s \circ \Phi^{-1}) \|_{m(s)} \\
				   & \leq & \|\phi_s\|_{m(s)+N} \\
				   &   =  & \|\Psi^{N,s+N} h\|_{m(s)+N} \\
				   & \leq &  (1 + t)^{2s+N}C''\bigl (\epsilon \cdot \|\Phi\|_{m(s) + 2N + n, 2} + \epsilon^{-(r + s + n/2)/(n/2)} \cdot  \|\Phi\|_{0,2} \bigr ) \|h\|_{s+2N} \, .   \\
\end{eqnarray*}
In the second inequality we use Condition $R_t$.  In the third inequality here we use Lemma 6.4.
That completes the argument.
\endpf 
\smallskip \\

\noindent {\bf Proof of Theorem 1.1:}  The last lemma tells us that 
$U \cdot (h \circ \Phi^{-1}) \in H^\infty (\overline{\Omega}_2)$ if $h \in H^\infty(\overline{\Omega}_1)$. 
Taking $h \equiv 1$, we conclude immediately that 
$$
\|U\|_{s, 2} \leq C''(1 + t)^s \bigl (1 + \epsilon \cdot \|\Phi\|_{m(s) + 2N + n, 2} + \epsilon^{-(m(s) + 2N + n/2)/(n/2)} \cdot \|\Phi\|_{0,2} \bigr ) \, .
$$

Next take $h = w_j$, where $w_j$ is the $j^{\rm th}$ coordinate function on $\Omega_2$.  We conclude now
that 
$$
\|U \cdot \bigl (\Phi^{-1}\bigr )_j \|_{s, 2} \leq C''(1 + t)^s \bigl (1 + \epsilon \cdot \|\Phi\|_{m(s) + 2N + n, 2} + \epsilon^{-(m(s) + 2N + n/2)/(n/2)} \cdot \|\Phi\|_{0,2} \bigr ) \, .   \eqno (6.7)
$$											      

Fix a point $z \in \Omega_1$.  The fact that $\Phi: \Omega_1 \ra \Omega_2$ is Lipschitz of order greater
than $(n-1)/n$ tells us that 
$$
|\nabla \Phi(z)| \leq C \cdot \delta_1(z)^{-1/n + \epsilon} 
$$
for some $\epsilon > 0$.  
Hence the complex Jacobian determinant $u$ of $\Phi$ at $z$
is bounded by $d_1^{-1 + \epsilon'}(z)$ for some 
$\epsilon' > 0$.  We know from results of Range [RAN1], proved with a direct
application of Hopf's lemma, that $d_1(\Phi^{-1}(w))^m \leq C \cdot d_2(w)$ for some positive integer $m$.
But in fact the bi-Lipschitz condition of order exceeding $(n-1)/n$ guarantees that $m$ must be 1.  

It follows then that $U$ must be of size
at least $d_2^{1 - \epsilon'}$ if it vanishes at some point of $\partial \Omega_2$ (it cannot
vanish in the interior).  That contradicts the smoothness of $U$ to the boundary.  
We conclude then that $U$ cannot vanish.  Hence it is bounded from 0 in modulus.
So we may see from (6.7) that 
\begin{align}
\|\bigl ( \Phi^{-1} \bigr )_j \|_{s, 2} & \leq  C''(1 + t)^s \biggl ( \biggl [ \frac{\epsilon_s}{(1 + 2t)^{N + n/2}} \biggr ]^{-0} \cdot \|\Phi\|_{0,2}  \notag \\
                                   &       \qquad + \frac{\epsilon_s}{(1 + 2t)^{N + n/2}} \cdot \|\Phi\|_{2s} \biggr ) \, .  \tag*{(6.8)}    \\	 \notag
\end{align}

Of course a similar argument may be applied with $\Phi^{-1}$ replace by $\Phi$ and the
roles of $\Omega_1$ and $\Omega_2$ reversed to see that 
\begin{align}
\|\bigl ( \Phi \bigr )_j \|_{s, 2} & \leq  C''(1 + t)^s \biggl ( \biggl [ \frac{\epsilon_s}{(1 + 2t)^{N + n/2}} \biggr ]^{-0} \cdot \|\Phi^{-1}\|_{0,2}  \notag \\
                                   &       \qquad + \frac{\epsilon_s}{(1 + 2t)^{N + n/2}} \cdot \|\Phi^{-1}\|_{2s} \biggr ) \, .  \tag*{(6.9)}    \\	 \notag
\end{align}

Now let $\lambda_\ell = 10^{- \ell}$.

In inequality (6.8), replace
$s$ by $\ell$ and multiply through by $\lambda_\ell$.   Likewise, in inequality (6.9), replace
$s$ by $\ell$ and multiply through by $\lambda_\ell$.  Now sum over $\ell$.  The result is
\begin{eqnarray*}
\lefteqn{\sum_\ell \biggl [ \lambda_\ell \|\bigl ( \Phi^{-1}\bigr )_j \|_{\ell, 2}  + \lambda_\ell \|\bigl ( \Phi \bigr )_j \|_{\ell, 2} \biggr ]}  \\
                & \leq & \sum_\ell \biggl [C'' \biggl ( 1 + \epsilon_\ell^{- 0} (\cdot  \|\Phi\|_{0,2} + \epsilon_\ell^{-0} \cdot  \|\Phi^{-1}\|_{0,2} \\
                &&  + \epsilon_\ell \|\Phi^{-1}\|_{2s,2} + \epsilon_\ell \|\Phi\|_{2s,2}\biggr ) \lambda_\ell \biggr ] \, . \\
\end{eqnarray*}
What is nice about this inequality is that we can now absorb the two $\epsilon_\ell$ terms
on the righthand side into the lefthand side.   In order to do this, we must note that
the term $\| \ \ \|_{2\ell, 2}$ on the righthand side has coefficient $\epsilon_\ell \lambda_\ell$
while the same term on the lefthand side has coefficient $\lambda_{2\ell}$.  So we must
choose $\epsilon_\ell$ so that $\epsilon_\ell \lambda_\ell \leq (1/2) \lambda_{2\ell}$.  Clearly
$\epsilon_\ell = (1/2) 10^{-\ell}$ will do the job.

The result is that
\begin{eqnarray*}
\lefteqn{\sum_\ell \biggl [ \lambda_\ell \|\bigl ( \Phi^{-1}\bigr )_j \|_{\ell, 2}  + \lambda_\ell \|\bigl ( \Phi \bigr )_j \|_{\ell, 2} \biggr ]}  \\
                & \leq & \sum_\ell C''' \lambda_\ell \cdot \bigl ( \epsilon_\ell^{-0} \cdot \|\Phi\|_{0,2} + \epsilon_\ell^{-0} \cdot \|\Phi^{-1}\|_{0,2} \bigr ) \, .  \\
\end{eqnarray*}

We may conclude from this last inequality that $\|\bigl ( \Phi^{-1}\bigr ))j\|_{j,2}$ and
$\|\bigl ( \Phi \bigr ))j\|_{j,2}$ are bounded.  Thus
the bihlomorphic mapping extends to a diffeomorphism of the closures.  That is our theorem.
\endpf
\smallskip \\
	     
We remark that, if we strengthen the hypotheses of our theorem
to $\Phi$ and $\Phi^{-1}$ both being Lipschitz 1, then it is immediate
that $u$ and $U$ are bounded and the proof simplifies notably.  Having
a Lipschitz condition of order less than 1 makes things a bit trickier.

\section{Concluding Remarks}

It would be natural to suppose that a theorem like the one that we prove here
is actually valid with only the assumption that $\Phi$ and $\Phi^{-1}$ are
Lipschitz of order $\epsilon$ for some $\epsilon > 0$.  The methods
that we have do not suffice to establish such a result.

We repeat that, of course, the hope is that no Lipschitz hypothesis should
be needed.  The conclusion should be true all the time.  That question will
be a topic for future research.

\newpage

\noindent {\Large \sc References}
\vspace*{.2in}

\begin{enumerate}

\item[{\bf [ADA]}]   R. A. Adams, {\it Sobolev Spaces}, Academic Press, New York, 1975.

\item[{\bf [BEL1]}]  S. R. Bell, Biholomorphic mappings and the $\dbar$ problem,
{\em Ann. Math.}, 114(1981), 103-113.

\item[{\bf [BEL2]}] S. R. Bell, Local boundary behavior of
proper holomorphic mappings, {\it Complex Analysis of Several
Variables} (Madison, Wis., 1982), 1--7, Proc.\ Sympos.\ Pure
Math., 41, Amer.\ Math.\ Soc., Providence, RI, 1984.

\item[{\bf [BELL]}] S. R. Bell and E. Ligocka, A simplification
and extension of Fefferman's theorem on biholomorphic
mappings, {\em Invent. Math.} 57(1980), 283--289.
	      
\item[{\bf [BUR]}] R. B. Burckel, {\it An Introduction to
Classical Complex Analysis}, Academic Press, New York, 1979.

\item[{\bf [DIF]}]  K. Diederich and J. E. Forn\ae ss, Pseudoconvex domains: 
Bounded strictly plurisubharmonic exhaustion functions, {\it Invent. Math.}
39(1977), 129--141. 

\item[{\bf [FEF]}] C. Fefferman, The Bergman kernel and
biholomorphic mappings of pseudoconvex domains, {\em Invent.
Math.} 26(1974), 1--65.
				
\item[{\bf [FOK]}] G. B. Folland and J. J. Kohn, {\it The Neumann Problem
for the Cauchy-Riemann Complex}, Princeton University Press, Princeton,
1972.						      

\item[{\bf [FRI]}] B. Fridman, Biholomorphic transformations that do not
extend continuously to the boundary, {\it Michigan Math.\ J.} 38(1991),
67--73.

\item[{\bf [GOL]}] G. M. Goluzin, {\it Geometric Theory of
Functions of a Complex Variable}, American Mathematical
Society, Providence, 1969.

\item[{\bf [GRK]}] R. E. Greene and S. G. Krantz, {\it Function
Theory of One Complex Variable}, 3rd ed., American
Mathematical Society, Providence, RI, 2006.

\item[{\bf [HEN]}] G. M. Henkin, An analytic polyhedron is not
holomorphically equivalent to a strictly pseudoconvex domain, (Russian)
{\it Dokl.\ Akad.\ Nauk SSSR} 210(1973), 1026--1029.

\item[{\bf [HOR]}]  L. H\"{o}rmander, $L^2$ estimates and existence theorems
for the $\dbar$ operator, {\it Acta Math.} 113(1965), 89--152.
							       w
\item[{\bf [KOH]}] J. J. Kohn, Global regularity for $\overline{\partial}$
on weakly pseudo-convex manifolds, {\it Trans.\ AMS} 181(1973), 273--292.

\item[{\bf [KRA1]}]  S. G. Krantz, {\it Function Theory of Several
Complex Variables}, 2nd ed., American Mathematical Society,
Providenc, RI, 2001.

\item[{\bf [RAN1]}] R. M. Range, The Carath\'{e}odory metric
and holomorphic maps on a class of weakly pseudoconvex
domains, {\it Pacific J. Math.} 78(1978), 173--189.

\item[{\bf [RAN2]}] R. M. Range, A remark on bounded strictly
plurisubharmonic exhaustion functions, {\it Proc.\ AMS}
81(1981), 220--222.

\item[{\bf [ROM]}] S. Roman, The formula of Fa\`{a} di Bruno, {\it Am.
Math. Monthly} 87(1980), 805-809.

\item[{\bf [STE]}] E. M. Stein, {\it Harmonic Analysis: Real
Variable Methods, Orthogonality, and Oscillatory Integrals},
Princeton University Press, Princeton, NJ, 1993.

\end{enumerate}
\vspace*{.25in}

\begin{quote}
Steven G. Krantz  \\
Department of Mathematics \\
Washington University in St. Louis \\
St.\ Louis, Missouri 63130  \\
{\tt sk@math.wustl.edu}
\end{quote}

\end{document}